\documentclass[3p,authoryear]{elsarticle}

\usepackage{lineno,hyperref}
\modulolinenumbers[5]

\journal{Mechanics Research Communications}

\usepackage[T1]{fontenc}
\usepackage[latin1]{inputenc}
\usepackage{amsmath}
\usepackage{mathrsfs}
\usepackage{amssymb}

\newtheorem{theorem}{Theorem}
\newtheorem{lemma}{Lemma}
\newtheorem{proposition}{Proposition}
\newtheorem{remark}{Remark}

\newcommand{\al}{\alpha}
\newcommand{\be}{\begin{equation}}
\newcommand{\ee}{\end{equation}}
\newcommand{\half}{\frac{1}{2}}
\newcommand{\p}{\partial}

\begin{document}

\begin{frontmatter}

\title{A thermoelastic theory with microtemperatures of type III}

\author[tun]{Moncef Aouadi\corref{cor1}}
\ead{moncefaouadi00@gmail.com}

\author[sal1]{Michele Ciarletta}
\ead{mciarletta@unisa.it}

\author[sal2]{Francesca Passarella}
\ead{fpassarella@unisa.it}

\cortext[cor1]{Corresponding author  at : Universit\'e de Carthage,
Ecole Nationale d'Ing\'enieurs de Bizerte, 7035, BP66, Tunisia.}

\address[tun]{Universit\'e de Carthage, UR Syst\`emes dynamiques et
applications, 17ES21\\ Ecole Nationale d'Ing\'enieurs de Bizerte,
7035, BP66, Tunisia}

\address[sal1]{Università degli Studi di Salerno, Dipartimento di Ingegneria Industriale\\
Via Giovanni Paolo II 132, 84084 Fisciano (SA), Italy}

\address[sal2]{Università degli Studi di Salerno, Dipartimento di Matematica\\
Via Giovanni Paolo II 132, 84084 Fisciano (SA), Italy}

\begin{abstract}
In this paper, we use the Green-Naghdi theory of thermomechanics of
continua to derive a nonlinear theory of thermoelasticity with
microtemperatures of type III.  This theory permits propagation of
both thermal and microtemperatures waves at finite speeds  with
dissipation of energy. The equations of the linear theory are also
obtained.   With the help of the semigroup theory of linear
operators  we establish that the linear anisotropic problem is well
posed and we study the asymptotic behavior of the solutions.
Finally,  we investigate the impossibility of the localization in
time of  solutions.

Published as:\\
Moncef Aouadi, Michele Ciarletta \& Francesca Passarella, 
\emph{Thermoelastic theory with microtemperatures and dissipative thermodynamics},
Journal of Thermal Stresses, \textbf{41}(4), 2018 522--542 \\
\url{https://doi.org/10.1080/01495739.2017.1383219}
\end{abstract}


\end{frontmatter}


\section{Introduction}
The origin of the theories of microtemperatures in elastic solids
goes back to works of \citet{Eringen1}, \citet{Grot},
\citet{Wozniak1, Wozniak2}, \citet{Iesan1, Iesan2}
and \citet{IQ1, IQ2}.   The experimental observations
have shown that the classical continuum theories cannot be used to
describe satisfactorily some phenomena. The interest of
microtemperatures is stimulated by the fact that this theory is
adequate to investigate important problems related to size effects
and nanotechnology. This represents an important improvement in the
perspective of real applications in many engineering and geophysics
applications and also in nanomanufactured activities.  Analogous to
the concept of micro-deformation introduced by 
\citet{Eringen1} in his celebrated theory of micromorphic continua,
\citet{Grot} extended the thermodynamics of a continuum with
microstructure  so that the points of a generic microelement are
assumed to have different temperatures. In \citep{Grot}, Grot
supposed that  the inverse of the microelement temperature is a
linear function of microcoordinates. The Clausius-Duhem inequality
is modified to include microtemperatures, and the first-order moment
of the energy equations are added to the usual balance laws of a
continuum with microtemperatures. The concept of microtemperatures
has been introduced for the first time by  \citet{Wozniak1,
Wozniak2}. He considered in \citep{Wozniak2} that the continuum
(macromedium)   is composed of particles $X$ with macrocoordinates.
Each particle is assumed at the same time as the origin of the
so-called system of spatial local microcoordinates $(x_1, x_2,
x_3)$. The spatial coordinates $x'_i$ of the point $X'$ of the
microelement  are represented in the form $x'_i = x_i
\psi_{ik}\xi_k$, where $x_i$ are the spatial coordinates of the
centroid $X$ of the microelement; $X'_k$ and $X_k$ are the material
coordinates of $X'$ and $X$, and $\xi_k = X'_k - X_k$. The functions
$\psi_{ik}$ are called microdeformations.  \citet{Wozniak1,
Wozniak2} assumed that all the properties of the macromedium of
particles $X$ are deduced from the structural properties of the
micromedium associated with each particle $X$ and from the way in
which these micromedia are interconnected. In particular, the
temperatures at points belonging to micromedium are assumed to be
known functions (depending on local microcoordinates) expressed by
the mean temperature of the micromedium and its gradient. In
\citep{Grot}, Grot established the thermodynamics of continua with
microstructure when the points of a generic microelement have
different temperatures. In this theory the temperature $\theta'$ at
the point $X'$ of the microelement  is a linear function of the
microcoordinates $x_k$, of the form $\theta'=\theta+ \tau_k\xi_k$,
where $\theta$ is the temperature at the centroid $X$. The vector
with the components $T_k$ defined by $T_k =-\tau_k/\theta$ is called
the microtemperatures vector.

The usual theory of heat conduction based on the classical Fourier's
law allows the phenomena of ``infinite diffusion velocity'' which is
not well accepted from a physical point of view. This  paradox of
the heat conduction, is physically unrealistic since it implies the
propagation of thermal waves with infinite speed.
 In contrast to the classical Fourier's law,
nonclassical thermal laws came into existence during the last
decades to eliminate this shortcomings. A survey article of
representative  these new models is due to \citet{[3]}. 
One of these theories may be mentioned that of Green and
Naghdi who have developed in a series of articles (see
\citep{GN1,GN2,GN3}) a thermomechanical theory of deformable continua
that relies on an entropy balance law rather than an entropy
inequality.  However, we want to mention the total compatibility of
the entropy balance law with the entropy inequality. They proposed
the use of the thermal displacement
\[\al({\bf x},t)=\int_{t_0}^{t}\theta({\bf x},s)ds+\al_0,\]
where $\theta$ is the empirical temperature, and considered three
theories labelled as type I, II and III, respectively. These
theories were based on an entropy balance law rather than the usual
entropy  inequality.
 The type I thermoelasticity coincides with the classical one; in
type II, the heat is allowed to propagate by means of thermal waves
but without dissipating energy and, for this reason, it is also
known as thermoelasticity without energy dissipation.   The heat
equation of type III, where the heat flux is a combination of type I
and II, contains both type I and II as limiting cases.
 In addition, the thermoelasticity of type III allows the constitutive functions
for free energy, stress tensor, entropy and heat flux to depend on
the strain tensor, the time derivative of the thermal displacement,
the  gradient  of  thermal  displacement  and  the  time  derivative
of the gradient of thermal displacement. This theory allows the
dissipation  energy,  but  the  heat  flux  is  partially determined
from  the  Helmholtz  free  energy  potential. Both, type II and
III, overcome the unnatural property  of  Fourier's law of infinite
propagation speed and  imply a finite wave propagation.

It is worth citing several papers on  thermoelastic theory with
microtemperatures and/or in the frame of Green-Naghdi models have
been published in the recent years (see, e.g., \citep{IQ1,
IQ2,Puri,QQ,Giorgi,[10],[11],[12],[13],[14],ACT,ACI}). One of these,
one can cite an intersting paper  published by \citet{IQ2} who
  derived a linear
theory of thermoelastic bodies with microstructure and
microtemperatures in the frame of Green-Naghdi theory of type II.
This theory permits propagation of both thermal and
microtemperatures waves at finite speeds; but without energy
dissipation and consequently with constant energy. It is worth
mentioning that  few applicability studies have been developed
 type II model. However, mathematical and physical analysis is needed to clarify
its applicability. In this paper we extend this theory to the
nonlinear type III model without considering the microstructure of
the material which is not our main concern here. The obtained
equations for temperature and microtemperatures allow the
transmission of heat with finite speeds and with energy dissipation.
For this, the type III model seems physically more acceptable than
type II model since it allows the propagation of the resulting
deformations with decaying energy when  the tensors of thermal and
microthermal conductivity are both positive definite. Also a
comparative test in \citep{Giorgi} between type III  and type I
models reveals that the type III model is more preferable . On the
other hand, we need to warm the reader that our model predicts an
asymptotic behavior of solution (see the fifth section) which is not
possible in the frame of type II model established in \citep{IQ2}.
 That is,  we prove that the solution  tends to zero when the time tends to infinity (see Theorem 3 in the fifth section).
Moreover, it is known that in many situations the decay can be
exponential. This proof is omitted in this paper for the sake of
brevity.  It should be noted that this paper is the first one  that
derives a Green-Naghdi theory of type III in the presence of
microtemperatures. Hence, the theory derived in this article will be
of great importance in real applications in many engineering and
geophysics applications and also in nanomanufactured activities.

The organization of this paper is as follows. In Section 2 we use
the theory established by \citet{GN1, GN2, GN3} to
obtain a nonlinear theory of thermoelastic with microtemperatures of
type  III, which admits the possibility of ``second sound'' with
energy dissipation. The type II model is also derived  as a step in
the derivation of the type III which our main concern in this paper.
In fact the type III model is obtained by adding to the heat flow of
the type II model the dependency on the temperature gradient and on
the microtemperatures gradient. The process of linearization of the
obtained equations is presented in Section 3.
 With the help of the semigroup
theory of linear operators an existence result is obtained in
Section 4. In Section 5, the asymptotic behavior for the solutions
of type III problem is studied. Finally,  in Section 6, we
investigate the impossibility of the localization in time of
solutions of type III problem. It is worth noting that we focus on
the analysis of the qualitative properties of solutions of type III
problem. However, some particular aspects of the type II problem are
also pointed out.

\section{Nonlinear  theory}
 In this section we present a nonlinear
theory of thermoelasticity with microtemperatures  in the context of
the Green-Naghdi model of type III.

We consider a continuous body  that at time $\tau_0$ occupies a
bounded region $\Omega$ of the Euclidean three-dimensional space
with smooth boundary   $\p \Omega$. For any sub-body   we denote
with $B$ the corresponding region in the reference configuration
$\Omega$, which is bounded by a regular surface  $\p {B}$, and with
${n_i}$ the components of the  unit outward normal  to $\p {B}$. We
shall employ the usual summation and differentiation conventions:
Latin subscripts are understood to range over the integers $(1, 2,
3)$, summation over repeated subscripts is implied and subscripts
preceded by a comma denote partial differentiation with respect to
the corresponding cartesian coordinate. Moreover a superposed dot
denotes the partial time derivative.

We take the configuration $\Omega$ as reference configuration and
refer the motion of the continuum to the reference configuration.
Fixed system of rectangular cartesian axes, we denote with $X_i$ the
coordinates of a point in the reference configuration and with $x_i$
the coordinates of the same point at time $t$, where
$x_i=\hat{x}_i(X_1, X_2, X_3, t)$. We assume that the functions
$\hat{x}_i$ are continuously differentiable as it is necessary.
 Following \citet{Grot} and \citet{Eringen1},
 we wish to extend the linear thermoelastic theory with microtemperatures of type II derived by 
  \citet{IQ2} to the nonlinear type III model. This can be done by deriving the balance laws  to the case in
which the temperatures  of the particles are different. In the
spirit of their works, let assume that $\textbf{X}$ be the center of
mass of a generic microelement in the reference configuration. We
assume that deriving the particle temperature $T'$ have the form
\citep{Grot}\be\label{temp0} {T'}(\textbf{X}',t)= {T}(\textbf{X},t) +
{{T}_i}(\textbf{X},t)({X}'_i-X_i)\ee   where the function ${T}_i$ is
called microtemperatures. The variables in \eqref{temp0} form a
vector and represent the variation of the temperature  within a
microvolume. They will be considered as independent thermodynamic
variables to be determined by the balance laws.

The balance of linear momentum can be written in the form \be
\label{motion}
  \rho  \ddot x_i = t_{ki,k}+\rho f_i \ee
where $t_{ik}$ is the first Piola Kirchhoff stress, $\rho$ is the
reference mass density and $f_i$ is the body force.  The surface
traction $t_i$ at regular points of $\p B$ are given by
\be\label{motion1}t_i= t_{ik}n_k.\ee

Following \citep{GN1, GN2, GN3} we postulate, for each microelement,
the local balance of entropy $$\rho' \dot S'= \Phi_{k,k}'+ \rho'
(s'+\xi')$$ where $\rho'$ is the mass density of the microelement,
$S'$ the entropy density per unit mass of the microelement,
$\Phi_k'$ the microentropy flux vector, $s'$ the external rate of
supply of entropy per unit mass of the microelement and $\xi'$ the
internal rate of production of entropy per unit mass of the
microelement. The heat flux vector associated with the microelement
is given by $q_i' =T'\Phi_i' $ and the external rate of supply of
heat per unit mass is defined by $r' = T's'$.

According to \citep{GN1, GN2, GN3}, we postulate the local balance of
entropy \be\label{localentropy}\rho \dot S= \Phi_{k,k}+ \rho
(s+\xi)\ee and the balance of first moment of entropy
\be\label{firstmoment} \rho\dot \varepsilon_i= \Lambda_{ji,j} +
\Phi_i - H_i + \rho (Q_i+\xi_i)\ee where $S$ is the entropy per unit
mass of the body, $\Phi_i$ is the entropy flux vector, $s$ is the
external rate of supply of entropy per unit mass, $\xi$ is the
internal rate of production of entropy per unit mass,
$\varepsilon_i$ is the first entropy moment vector, $\Lambda_{ij}$
is the first entropy flux moment tensor, $H_i$ is the mean entropy
flux vector, $Q_i$ is the first moment of the external rate of
supply of entropy and $\xi_i$ is the first moment of the internal
rate of production of entropy.

Moreover, following the arguments of \citet{GN1} and
\citet{IQ2}, we postulate for every
subregion $B$ of $\Omega$ and every time $t$, the following  balance
equations for the energy
 \be
 \label{energy}
 \int_B \rho (\dot x_i
\ddot x_i+ \dot e)dv=\int_B \rho (f_i \dot x_i
+sT+Q_iT_i)dv+\int_{\p B} (t_i\dot x_i+\Phi T+\sigma_iT_i) da, \ee
where $e$ is the internal energy per unit mass.  The entropy flux
$\Phi$ and the first entropy moment flux vector $\sigma_i$ at
regular points of $\p B$ are given by \be \label{T_q_Phi_lambda0}
{\Phi}=\Phi_k{n}_k\,,\qquad \sigma_k= \Lambda_{jk}n_j.\ee Therefore,
thanks to the arbitrariness of $B$, equations
\eqref{motion}-\eqref{energy},  we obtain the following local form
for the balance equations \be\label {local_form}
 \rho \dot e=t_{ki} \dot x_{i,k}+\rho s T+\rho Q_iT_i+\Phi_{k,k}T+\Phi_k T_{,k}+\Lambda_{kj,k}T_{j}+\Lambda_{kj}T_{j,k}.
\ee By eliminating $\rho s T$  and $\Lambda_{kj,k}T_{j}$ from
\eqref{local_form} through the use of $(\ref{localentropy})$ and
(\ref{firstmoment}),  we obtain\be\label {local_form1}
 \rho \dot e=t_{ki} \dot x_{i,k}+\rho\dot ST-\rho\xi T+\Phi_k T_{,k}+\Lambda_{kj}T_{j,k}
 +\rho\dot  \varepsilon_kT_k+(H_k-\Phi_k)T_k-\rho\xi_kT_k.
\ee If we introduce the specific Helmholtz free energy per unit mass
\be\label{Psi} \Psi=e-TS-\varepsilon_kT_k, \ee
 the energy equation \eqref{local_form1} becomes
\be\label{Gibbs} \rho(\dot \Psi+\dot T S+\dot T_k
\varepsilon_k)=t_{ki} \dot x_{i,k}-\rho\xi T+\Phi_k
T_{,k}+\Lambda_{kj}T_{j,k}
 +(H_k-\Phi_k)T_k-\rho\xi_kT_k.
\ee We introduce the notations \be\label {Lebon}q_j = \Phi_j
T,\qquad q_{ji} = \Lambda_{ji}T. \ee The heat flux $q$ and the heat
flux moment vector $\Lambda_i$ at regular points of $\p B$ are given
by \be\label {Lebon1}q =q_jn_j, ,\qquad \Lambda_i =q_{ ji}n_j \ee
respectively. We shall now introduce the nonlinear model of type II
to deduce the type III model.
\subsection{Type II - dissipationless theory}According to the Green-Nagdhi theory \citep{GN1}, we introduce
the thermal displacement $\al$ whose derivative coincides with the
absolute temperature, i.e., $ \dot\alpha=T$. This scalar, on the
macroscopic scale, is regarded as representing some ``mean'' thermal
displacement magnitude on the molecular scale. In a similar way, we
introduce a  scalar function $\beta_i$ related to the
microtemperatures by the equation
  $\dot \beta_i=T_i$ (see \citep{GN1, IQ2}).

We assume that the response functions \be\label{response} \Psi,\
t_{kj},\ S,\ \varepsilon_i,\ \Phi_k,\ \Lambda_{jk},\ H_k,\ \xi ,\
\xi_k \ee depend on the set of the independent variables
\[\mathcal{A}_1=(x_{i,k},\ T,\ T_k,\ \al_{,k},\ \beta_{i,k}).\]
Thus, we consider constitutive equations of the form \be
\mathcal{F}=\widehat{ \mathcal{F}}(\mathcal{A}_1) \ee and we assume
that the response functions are of $C^1-$class.
Using the chain rule \be\label{chain}\dot \Psi=
 \frac{\p \widehat{\Psi}}{\p x_{j,k}}\dot x_{j,k} +
 \frac{\p\widehat{\Psi}}{\p T}\dot T+
\frac{\p \widehat{\Psi}}{\p T_k}\dot T_k+ \frac{\p\widehat{\Psi}}{\p
\alpha_{,k}}\dot\alpha_{,k}+ \frac{\p \widehat{\Psi}}{\p
\beta_{j,k}}\dot \beta_{j,k}, \ee the comparison of Eqs.
(\ref{Gibbs}) and (\ref{chain}) yields\begin{align}\begin{split}
\label{restrictions} &\left(\rho\frac{\p\widehat{\Psi}}{\p T}+\rho
\widehat S\right)\dot T+ \left(\rho\frac{\p\widehat{\Psi}}{\p
T_i}+\rho \widehat \varepsilon_i\right)\dot T_i+ \left(\rho\frac{\p
\widehat{\Psi}}{\p x_{j,k}}-\widehat t_{kj}\right)\dot x_{j,k} +
\left(\rho\frac{\p\widehat{\Psi}}{\p \alpha_{,k}}-\widehat
\Phi_{k}\right)\dot\alpha_{,k}\cr&+ \left(\rho\frac{\p
\widehat{\Psi}}{\p \beta_{i,k}}-\widehat \Lambda_{ki}\right)\dot
\beta_{i,k}+  \rho T \widehat \xi+(\widehat \Phi_k-\widehat
H_k)T_k+\rho\widehat \xi_kT_k=0
\end{split} \end{align}
 which must hold for all choice of $\dot T$, $\dot T_i$, $\dot x_{j,k}$, $\dot\alpha_{,k}$ and $\dot \beta_{i,k}$.
 From this
equality we see that the constitutive equations  are compatible with
the energy equation if satisfy the following relations
 \be\label{restrictions1}
 \begin{array}{l}
\displaystyle{ \ \ \ \Psi=\widehat{\Psi}(\mathcal{A}_1),\ \ \
S=-\frac{\p\widehat{\Psi}(\mathcal{A}_1)}{\p T}, \ \ \
\varepsilon_i=-\frac{\p\widehat{\Psi}(\mathcal{A}_1)}{\p T_i}, \ \ \
t_{kj}=\rho\frac{\p \widehat{\Psi}(\mathcal{A}_1)}{\p x_{j,k}},
 }\\
\\
\displaystyle{\Phi_{k}= \rho\frac{\p
\widehat{\Psi}(\mathcal{A}_1)}{\p \alpha_{,k}},\ \ \
\Lambda_{ki}=\rho\frac{\p
\widehat{\Psi}(\mathcal{A}_1)}{\p\beta_{i,k}},\ \ \   \rho T
 \xi+( \Phi_k- H_k)T_k+\rho
\xi_kT_k=0.}
  \end{array}
\ee \begin{remark}\label{rq1} Since the last equation of
\eqref{restrictions1} must be satisfied for all process and the
temperature and the microtemperatures can not vanish for all
process, we conclude that
$${\xi}=0\qquad \textrm{whenever}\qquad   \Phi_k- H_k+\rho
\xi_k=0 .$$\end{remark}
 The thermal displacement $\al$ and the microtemperatures
displacement $\beta_k$ are defined analogously to the well-known
mechanical displacement. For this,
 the entropy flux vector $\Phi_k$  and the first entropy flux moment tensor $\Lambda_{ij}$
are deduced from a potential  in the same way as the stress tensor
is derived in mechanics.
\subsection{Type III - dissipation theory}
Whereas in case of heat flow of type II the response functions
(\ref{response}) are assumed to depend on the { material deformation
gradient} $x_{i,k}$, the temperature  $T$,  the microtemperatures
$T_i$, the thermal displacement gradient $\al_{,k}$ and the
microtemperatures displacement gradient $\beta_{i,k}$, for type III
we now add the dependency on the temperature gradient $\dot\al_{,k}$
and on the microtemperatures gradient $\dot\beta_{i,k}$. Hence, we
assume that the response functions (\ref{response}) depend on  the
set of the independent variables
\[\mathcal{A}_2=(x_{i,k}, T, T_i, \al_{,k},\beta_{i,k}, \dot\al_{,k},\dot\beta_{i,k})=(x_{i,k}, T, T_i,
 \al_{,k},\beta_{i,k}, T_{,k},T_{i,k}).\]
In this case, using the chain rule \be\label{chain1}\dot \Psi=
 \frac{\p \widehat{\Psi}}{\p x_{j,k}}\dot x_{j,k} +
 \frac{\p\widehat{\Psi}}{\p T}\dot T+
\frac{\p \widehat{\Psi}}{\p T_k}\dot T_k+ \frac{\p\widehat{\Psi}}{\p
\alpha_{,k}}T_{,k}+ \frac{\p \widehat{\Psi}}{\p \beta_{i,k}}T_{i,k}+
\frac{\p\widehat{\Psi}}{\p T_{,k}}\dot T_{,k}+ \frac{\p
\widehat{\Psi}}{\p T_{i,k}}\dot T_{i,k}, \ee the comparison of Eqs.
(\ref{Gibbs}) and (\ref{chain1}) yields
\[\frac{\p\widehat{\Psi}}{\p T_{,k}}=0\,,\quad \frac{\p\widehat{\Psi}}{\p T_{i,k}}= 0,
\]
that is $\Psi=\widehat{\Psi}(x_{i,k}, T, T_k,
\al_{,k},\beta_{i,k})=\widehat{\Psi}(\mathcal{A}_1)$,  and
\[\ \ \ S=-\frac{\p\widehat{\Psi}(\mathcal{A}_1)}{\p T}, \ \ \
\varepsilon_i=-\frac{\p\widehat{\Psi}(\mathcal{A}_1)}{\p T_i}, \ \ \
t_{kj}=\rho\frac{\p \widehat{\Psi}(\mathcal{A}_1)}{\p x_{j,k}}, \]
\be\label{restrictions2} \ \ \ \left(\rho\frac{\p \widehat{\Psi}}{\p
\alpha_{,k}}- \widehat{\Phi}_{k}\right) T_{,k}+\left(\rho\frac{\p
\widehat{\Psi}}{\p\beta_{i,k}}- \widehat{\Lambda}_{ki}\right)
T_{i,k}+ \rho T
 \xi+( \Phi_k- H_k)T_k+\rho
\xi_kT_k=0. \ee

\section{ Linear theory}

 We consider a reference configuration  which is
 in thermal  equilibrium and free from stresses, with
 $\alpha$ and $\beta_k$ constant. We assume that the deformations and the changes of temperature and
 microtemperatures
are very small with respect to the reference configuration in such
way that, if $T_0$  and $T_i^0$ are respectively the (constant)
absolute temperature and the (constant) absolute microtemperatures
of the body in the reference configuration,
 we can write
  \be x_i-X_i= u_i =\varepsilon u'_i,   \ \ \
\ T-T_0=\theta=\varepsilon\theta',\ \ \ T_i-T_i^0= M_i=\varepsilon
M_i' \ee where $\varepsilon$ is a constant small enough for squares
and higher powers to be neglected, and $ u'_i$, $\theta'$ and  $
M_i'$ are independent on $\varepsilon$. Under these hypotheses, the
strain tensor is approximated with
\be\label{geometric}e_{ik}=\frac12(u_{i,k}+u_{k,i}).\ee
\subsection{Green-Naghdi theory of type II}
The   set of the independent variables for the Green-Naghdi model of
type II (without energy dissipation) becomes
\[\widetilde{\mathcal{A}}_1=(e_{ik}, \theta, M_i,\tau_{,i},R_{i,k})\]
where
\[ \tau=\int_{t_0}^t \theta ds,\ \ \
\ R_i=\int_{t_0}^t M_i ds.\] To obtain a linear theory, we consider
the free energy $\Psi$ function in the quadratic approximation
\be\label{Psi_lin} \rho\Psi=\half A_{ijkl} e_{ij}e_{kl}-a_{ij}
e_{ij}\theta+ B_{ijkl} e_{ij}R_{k,l}-b_{ij} R_{i,j}\theta
-\frac{\rho c_E}{2T_0}\theta^2+\half K_{ij}\tau_{,i}\tau_{,j}+\half
C_{ijkl} R_{i,j}R_{k,l} -d_{ij}M_{i}\tau_{,j}-\half c_{ij}M_{i}M_{j}
\ee where  $c_E$ is the specific heat at constant strain, $A_{ijkl}$
is the tensor of elastic constants, $a_{ij}$ is the tensor of
thermal expansion, $B_{ijkl}$ is the tensor of microtemperatures
expansion,  $K_{ij}$ is the tensor of thermal conductivity and
$c_{ij} $ is the tensor of microthermal conductivity. The tensors
$d_{ij}$ and $b_{ij} $ are, respectively, measure of thermal and
microthermal gradient displacement. The constitutive coefficients
have the following symmetries \be\label {sym1}
\begin{array}{lll}A_{ijkl} &=& A_{jikl} =
A_{kli j},\ a{_{i j}} = a_{ ji} ,\ B_{ijkl} = B_{jikl} = B_{kli j},\
\ K_{i j} =K_{ ji} \cr C_{ijkl} &=& C_{jikl} = C_{kli j},\ c_{i j} =
c_{ ji} ,\ b_{i j} = b_{ ji} ,\ d_{i j} = d_{ ji}.\\
\end{array}
\ee

For centrosymmetric materials, the compatibility conditions
(\ref{restrictions1}) give the following linear constitutive
equations \be\label {linear1}
\begin{array}{lll}
t_{ij}&=&A_{ijkl}e_{kl}-a_{ij}\theta+B_{ijkl}R_{k,l}, \\
\rho  S&=&\displaystyle{a_{ij}e_{ij}+\frac{\rho c_E}{T_0}\theta+b_{ij}R_{i,j}},\\
\rho\varepsilon_i&=&c_{ij}M_{j}+d_{ij}\tau_{,j},\\
 \Phi_i&=&-d_{ij}M_{j}+K_{ij}\tau_{,j},\\
 \Lambda_{ij}&=&B_{klji}e_{kl}-b_{ji}\theta+C_{jikl}R_{k,l}.\\
\end{array}
\ee Moreover, by using (\ref{linear1}), the linear approximation of
(\ref{Lebon}) is given by
\begin{align}\begin{split}\label{Phi_GN2} q_i&= T_0
\Phi_i= T_0(-d_{ij}M_{j}+K_{ij}\tau_{,j}),\\
q_{ij}&= T_0 \Lambda_{ij}=
T_0(B_{klij}e_{kl}-b_{ji}\theta+C_{jikl}R_{k,l}).\end{split}\end{align}
By Remark \ref{rq1}, we find that, in the linear theory of
centrosymmetric materials, we have $\xi = 0$ and consequently $
\Phi_k- H_k+\rho \xi_k=0$. Thus, \eqref{localentropy} and
\eqref{firstmoment} become\be\label{newlocalentropy}\rho \dot S=
\Phi_{k,k}+ \rho s,\qquad \rho\dot \varepsilon_i= \Lambda_{ji,j}+
\rho Q_i.\ee

The basic equations of the linear theory of type II consist of the
equations of motion \eqref{motion} and \eqref{motion1}, the energy
equations \eqref{newlocalentropy}, the constitutive equations
\eqref{linear1} and the geometrical equation \eqref{geometric}.
These equations furnish the following system of partial differential
equations for the unknown functions $u_i,\ \tau,\ R_i$
\begin{align}\begin{split}\label{linear0}
 \rho \ddot u_i &=(A_{jikl}e_{kl}-a_{ji}\dot \tau+B_{jikl}R_{l,k})_{,j}+\rho f_i,
\\
c\ddot \tau&=-a_{ij}\dot e_{ij}+(-d_{ij}\dot
R_{j}+K_{ij}\tau_{,j})_{,i} -b_{ij}\dot R_{i,j} +\rho s,
\\
c_{ij}\ddot R_{j}&= (B_{klji}e_{kl}-b_{ji}\dot
\tau+C_{jikl}R_{k,l})_{,i}-d_{ij}\dot \tau_{,j}+\rho Q_i,
\end{split}\end{align}
where $c= \frac{\rho c_E}{T_0}$. If the material is isotropic, then
the constitutive equations \eqref{linear1} become
\be
\begin{array}{lll}\label{iso22}
 t_{ij} &=&\lambda e_{kk}\delta_{ij}+2\mu e_{ij}-\beta\theta\delta_{ij}+\gamma_1R_{k,k}\delta_{ij}+2\gamma_2(R_{i,j}+R_{j,i}),\\
 \rho  S&=&\displaystyle{\beta e_{kk}+\frac{\rho c_E}{T_0}\theta+\varpi  R_{k,k}},\\
 \rho  \varepsilon_i&=&\alpha M_{i}+\hbar \tau_{,i}\\
   \Phi_i&=& -\hbar M_{i}+K \tau_{,i},   \\
  \Lambda_{ij}&=&\gamma_1e_{kk}\delta_{ij}+2\gamma_2 e_{ij}-\varpi \theta\delta_{ij}+\eta_1R_{k,k}\delta_{ij}+
  \eta_2R_{i,j}+\eta_3R_{j,i},
\end{array}
\ee where $\delta_{ij}$ is the Kronecker's delta, $\lambda$ and
$\mu$ are Lam\'e's constants, $\beta=(3\lambda+2\mu)\alpha_t$ and
$\alpha_t$ is the coefficient of linear thermal expansion.

By following the previous procedure, it follows from \eqref{iso22}
that the field equations of the theory of homogeneous and isotropic
bodies for the functions $u_i,\ \tau,\ R_i$ can be expressed as
\begin{align}\begin{split}
 \rho \ddot u_i&= \mu\Delta u_i+
(\mu+\lambda) u_{j,ji} -\beta\dot \tau_{,i}+\gamma_2\Delta
R_i+(\gamma_1+\gamma_2)R_{j,ji}+ \rho f_i\cr c\ddot \tau&=-\beta\dot
 u_{r,r}+K\Delta\tau-(\varpi+\hbar)\dot R_{i,i}+\rho s\cr
  \alpha \ddot R_i&=\gamma_2\Delta u_i+
(\gamma_1+\gamma_2) u_{j,ji}+\eta_2\Delta R_i+ (\eta_1+\eta_3)
R_{j,ji}-\hbar \dot \tau_{,i} +\rho
Q_i,\nonumber\end{split}\end{align} where $\Delta$ is the Laplacian.

\begin{remark} System \eqref {linear0}  can be deduced from system (3.18) of  \citep{IQ2} by eliminating the
porosity.  The model of the Green-Naghdi theory of type III was not
covered in \citep{IQ2}. We derive this model here and study some of
the corresponding  qualitative properties to fill this
gap.\end{remark}
\subsection{ Green-Naghdi theory  of type III}

By using the same  hypotheses and notations of the previous
subsection, the set of the independent variables for the
Green-Naghdi model of type III  (with energy dissipation) becomes
\[\widetilde{\mathcal{A}}_2=(e_{ik}, \theta, M_i,\tau_{,k},R_{i,k}, \theta_{,k}, M_{i,k})\]
 and, as a consequence of (\ref{restrictions2}), the quadratic approximation of the free energy $\Psi$
 is given by (\ref{Psi_lin}).
 Moreover, the constitutive equations for $t_{kj}$, $S$ and $\rho\varepsilon_i$
  are equal to (\ref{linear1})$_{1,2,3}$, respectively. The condition (\ref{restrictions2}) leads, in the linear context, to
\be\label{res_III} \left(\rho\frac{\p \widehat{\Psi}}{\p \tau_{,i}}-
\widehat{\Phi}_{i}\right) \theta_{,i}+\left(\rho\frac{\p
\widehat{\Psi}}{\p R_{i,j}}- \widehat{\Lambda}_{ji}\right) M_{i,j}+
\rho T_0
 \xi+( \Phi_i- H_i)T^0_i+\rho
\xi_iT^0_i=0. \ee The condition (\ref{res_III}) is satisfied if we
choose for a centrosymmetric material
 \be
 \label{linear3}
 \begin{array}{lll}
  \Phi_i&=&-d_{ji}M_{j}+K_{ij}\tau_{,j}+\tilde{K}_{ij}\theta_{,j},\\
 \Lambda_{ij}&=&B_{klji}e_{kl}-b_{ji}\theta+C_{jikl}R_{k,l}+\tilde{C}_{jikl}M_{k,l},
 \\
q_i&=& T_0 \Phi_i
=T_0\left[-d_{ji}M_{j}+K_{ij}\tau_{,j}+\tilde{K}_{ij}\theta_{,j}\right],
\\
q_{ij}&=& T_0 \Lambda_{ij}
=T_0\left[B_{klji}e_{kl}-b_{ji}\theta+C_{jikl}R_{k,l}+\tilde{C}_{jikl}M_{l,k}\right],\\
 \rho T_0
 \xi&=&
\tilde{K}_{ij}\theta_{,i}\theta_{,j}+\tilde{C}_{jikl}M_{l,k}M_{i,j},
\\
\rho \xi_i&=& H_i- \Phi_i,
 \end{array}
 \ee
where $\tilde{K}_{ij}$  and $\tilde{C}_{ijkl}$  are tensors
characteristic of the type III model. In view of
$\eqref{linear3}_{5,6}$, the local balance of entropy  and the
balance of first moment of entropy, in the linear context, are given
by \eqref{newlocalentropy}. Using \eqref{linear3} and
\eqref{newlocalentropy}, we obtain the following evolutive equations
of the theory of thermoelastic centrosymmetric materials with
microtemperatures  of type III (with energy dissipation)

\begin{align}\begin{split}\label{linear01}
 \rho \ddot u_i &=(A_{jikl}e_{kl}-a_{ji}\dot \tau+B_{jikl}R_{l,k})_{,j}+\rho f_i,
\\
 c\ddot \tau&=-a_{ij}\dot e_{ij}+(-d_{ij}\dot R_{j}+K_{ij}\tau_{,j}+\tilde{K}_{ij}\dot\tau_{,j})_{,i} -b_{ij}\dot R_{i,j} +\rho s,
\\
c_{ij}\ddot R_{j}&= (B_{klji}e_{kl}-b_{ji}\dot
\tau+C_{jikl}R_{k,l}+\tilde{C}_{jikl}\dot R_{k,l})_{,j}-d_{ij}\dot
\tau_{,j}+\rho Q_i,
\end{split}\end{align}where the above constitutive coefficients
 satisfy the following symmetry
relations
 \be\label {sym}
\begin{array}{lll}A_{ijkl} &=& A_{jikl} =
A_{kli j},\ a{_{i j}} = a_{ ji} ,\ B_{ijkl} = B_{jikl} = B_{kli j},\
\ K_{i j} =K_{ ji} \cr C_{ijkl} &=& C_{jikl} = C_{kli j},\ c_{i j} =
c_{ ji} ,\ b_{i j} = b_{ ji} ,\ d_{i j} = d_{ ji},\ \ \
\tilde{K}_{ij}=\tilde{K}_{ji},\ \ \ \
\tilde{C}_{ijkl}=\tilde{C}_{klij}.
\end{array}
\ee

  Remark that the evolutive equations (\ref{linear0}) of the
thermoelastic diffusion theory of type II (without energy
dissipation) can be deduced from  (\ref{linear01}) by taking
$\tilde{K}_{ij}=\tilde{C}_{ijkl}=0$.

If the material is isotropic, then the constitutive equations become
\[
\begin{array}{lll}
 t_{ij} &=&\lambda e_{kk}\delta_{ij}+2\mu e_{ij}-\beta\theta\delta_{ij}+\gamma_1R_{k,k}\delta_{ij}+2\gamma_2(R_{i,j}+R_{j,i}),\\
 \rho  S&=&\displaystyle{\beta e_{kk}+\frac{\rho c_E}{T_0}\theta+\varpi  R_{k,k}},\\
 \rho  \varepsilon_i&=&\alpha M_{i}+\hbar \tau_{,i},\\
   \Phi_i&=& -\hbar M_{i}+K \tau_{,i}+H\theta_{,j},   \\
  \Lambda_{ij}&=&\gamma_1e_{kk}\delta_{ij}+2\gamma_2 e_{ij}-\varpi \theta\delta_{ij}+\eta_1R_{k,k}\delta_{ij}+
  \eta_2R_{i,j}+\eta_3R_{j,i}+\varrho_1M_{k,k}\delta_{ij}+
  \varrho_2M_{i,j}+\varrho_3M_{j,i}.
\end{array}
\]

By following the previous procedure, it follows that the field
equations of the theory of type III for homogeneous and isotropic
bodies for the functions $u_i,\ \tau,\ R_i$ are
\begin{align}\begin{split}\label{linear3b}
 \rho \ddot u_i&= \mu\Delta u_i+
(\mu+\lambda) u_{j,ji} -\beta\dot \tau_{,i}+\gamma_2\Delta
R_i+(\gamma_1+\gamma_2)R_{j,ji}+ \rho f_i,\cr c\ddot
\tau&=-\beta\dot
 u_{r,r}+K\Delta\tau+H\Delta\dot \tau-(\varpi+\hbar)\dot R_{i,i}+\rho s,\cr
  \alpha \ddot R_i&=\gamma_2\Delta u_i+
(\gamma_1+\gamma_2) u_{j,ji}+\eta_2\Delta R_i+ (\eta_1+\eta_3)
R_{j,ji}+\varrho_2\Delta\dot  R_i+ (\varrho_1+\varrho_3) \dot
R_{j,ji}-\hbar\dot \tau_{,i} +\rho Q_i.\end{split}\end{align}

 To the
field of equations \eqref{linear01} (or \eqref{linear0}) we add
initial and  boundary conditions. Summarizing, the following initial
boundary value problems are to be solved for type III model (or Type
II  model):

Find $(u_i,v_i, \tau,\theta,M_i,R_i)$ solution to \eqref{linear01}
(or \eqref{linear0}) subject to the initial conditions \be
\label{IC1}
\begin{array}{lll} u_i(\cdot,0)&=&u_i^0,\ \ \ v_i(\cdot,0)= v_i^0,\ \ \
 \tau(\cdot,0)=\tau^0,\\ \theta(\cdot,0)&=& \theta^0, \  \ \ M_i(\cdot,0)=M_i^0,\
 \ \ R_i(\cdot,0)= R_i^0  \ \ \hbox{in}\   \Omega, \end{array}
\ee and the boundary conditions \be\label{BC1}
\begin{array}{lll} u_i &=&\tilde{u_i}\ \hbox{on}\ \p \Omega_u \times
 (0,\infty),\ \ \ \ \ \tau =\tilde{\tau}\ \hbox{on}\ \p \Omega_\tau \times
 (0,\infty),\ \ \ \ \ \ M_i=\tilde{M}_i\ \hbox{on}\ \p \Omega_{M_i} \times
 (0,\infty),\\
t_{ji}n_j &=&\tilde{t_i}\ \hbox{on}\  \p \Omega_t \times
 (0,\infty),\ \ \   \Phi_i n_i =\tilde{\Phi}\ \hbox{on}\  \p \Omega_\Phi \times
 (0,\infty),\ \ \  \Lambda_{ji} n_j =\tilde{\sigma}_i\ \hbox{on}\  \p \Omega_{\sigma_i} \times
 (0,\infty),\end{array}
\ee where $\tilde{u_i},\ \tilde{\tau},\ \tilde{M}_i,\ \tilde{t_i},\
\tilde{\Phi}$ and $\tilde{\sigma}_i$ are prescribed functions,
$u_i^0,\ v_i^0,\  \tau^0,\  \theta^0,\  M_i^0$ and
 $R_i^0$ are given and
 \[\p \Omega=\p \Omega_u\cup\p \Omega_t=\p \Omega_{\tau}\cup \p \Omega_\Phi=\p \Omega_{M_i}\cup \p \Omega_{\sigma_i},\
  \ \hbox{and}\  \ \p \Omega_u\cap\p \Omega_t=\p \Omega_{\tau}\cap \p \Omega_\Phi=\p \Omega_{M_i}\cap \p \Omega_{\sigma_i}=
 \emptyset \,. \]

In the following,  some qualitative properties of  the solutions to
type III problem are studied.
 However, some particular aspects of the type II problem are also pointed out.

\section{ Well-posedness}

We will prove the existence, uniqueness and continuous dependence
from the initial values and the external loads of the solution for
system  (\ref{linear01}) using the semigroups theory. Seeking for
simplicity, we will restrict ourselves to homogeneous boundary
conditions \be \label{BC}{\bf{u}}=0,\ \ \ \ \tau=0,\ \ \ \
\textbf{R}=0,\quad  \hbox{on}\ \ \p \Omega \times (0,\infty)\ee
where $ \textbf{u}$ and $ \textbf{R}$ denote the vectors of
components $u_i$ and $R_i$, respectively.

In the rest of the paper we assume:
\begin{itemize}
  \item [ (i)] relations (\ref{sym}) are satisfied;
  \item [ (ii)] $\rho>0$ and\footnote{  The inequality sign is a consequence of the Second Law of Thermodynamics,
   which requires the non-negativeness of the functional $\xi$ (see \citep{GN1}) and of our choice (\ref{linear3})$_5$.}
\be\label{tilda} \tilde{K}_{ij}\theta_{,i}\theta_{,j}+
\tilde{C}_{ijkl}M_{l,k}M_{i,j}\geq 0\,; \ee
\item [ (iii)] there exists a positive constant $c_0$ such that
\be\label{positive} \int_\Omega\Big( A_{ijkl}e_{kl}
e_{ij}+2B_{ijkl}e_{ij}R_{k,l}+K_{ij}
\tau_{,i}\tau_{,j}+C_{ijkl}R_{i,j}R_{k,l} \Big)dv\ge
c_0\int_\Omega(e_{ij}e_{ij}+\tau_{,i}\tau_{,i}+R_{i,j}R_{i,j})dv.\ee
\end{itemize}
We now wish to transform the  initial boundary value  problem
defined by the equations (\ref{linear01}), the initial conditions
(\ref{IC1}) and the boundary conditions (\ref{BC})  to an abstract
problem on a suitable Hilbert space. In what follows we use the
notations $v_i=\dot {u_i},\ \theta=\dot \tau,\ M_i=\dot R_i$. Let
define the following Hilbert space:
 \[\mathscr{H}=\left\{(u_i, v_i, \tau,\theta, R_i,M_i);\
  u_i,\ R_i\in {\textbf{{W}}}_0^{1,2}(\Omega);\ v_i,\ M_i \in \textbf{{L}}^{2}(\Omega)
 ;\  \tau\in {{W}}_0^{1,2}(\Omega);\ \theta \in
 L^{2}(\Omega)\right\},\]where ${W}_0^{1,2}(\Omega)$ and $\ L^{2}(\Omega)$ are
 the familiar Sobolev spaces and
 \[{\textbf{W}}_0^{1,2}(\Omega)=[{{W}}_0^{1,2}(\Omega)]^3,\ \ \ \ \ {\textbf{L}}^{2}(\Omega)=[{L}^{2}(\Omega)]^3.\]
  We
introduce the inner product in $\mathscr{H}$ defined by
\be\label{inner} \langle \mathcal{U},\mathcal{U}^*\rangle =
\half\int_\Omega \Big(\rho {v_i}{v_i}^* +c \theta \theta^*+c_{ij}
M_i M_j^* +2\mathcal{W}[(u_i, \tau,R_i),(u_i^*,\tau^*,R_i^*)]
\Big)dv, \ee where $\mathcal{U}=(u_i, v_i, \tau,\theta, R_i,M_i)$,
$\mathcal{U}^*=(u_i^*, v_i^*, \tau^*,\theta^*, R_i^*,M_i^*)$ and
\[
2\mathcal{W}[(u_i, \tau,R_i),(u_i^*,\tau^*,R_i^*)]= A_{ijkl}e_{kl}
e^*_{ij}+B_{ijkl}(e_{ij}R^*_{k,l}+e^*_{ij}R_{k,l})+K_{ij}
\tau_{,i}\tau^*_{,j}+C_{ijkl}R_{i,j}R^*_{k,l}.
\]
If we recall the assumption (\ref{positive}), the first Korn
inequality and the Poincar\'e inequality, we conclude that
\[\int_\Omega{\mathcal{W}}[(\textbf{u},\tau,\textbf{R}),(\textbf{u}, \tau,\textbf{R})] dv\] defines
a norm that is equivalent to the usual norm in
${\textbf{W}}_0^{1,2}(\Omega)\times {{W}}_0^{1,2}(\Omega)\times
{\textbf{W}}_0^{1,2}(\Omega)$. Hence, the bilinear form
(\ref{inner}) defines an inner product equivalent to the usual one
in $\mathscr{H}$.

We introduce the following operators \be\label {matrix1}
\begin{array}{lll}  A_i \textbf{u}  &=& \rho^{-1}(A_{jikl}u_{k,l})_{,j},\\
B_i\theta &=& -\rho^{-1}(a_{ji}\theta)_{,j},\\
C_i\textbf{R} &=& \rho^{-1}(B_{jikl}R_{l,k})_{,j},\\
 D\textbf{v}&=&-c^{-1} a_{ij}v_{i,j},\\
 E\tau&=&c^{-1} (K_{ij}\tau_{,j})_{,i},\\
 G\theta&=&c^{-1} (\tilde{K}_{ij}\theta_{,j})_{,i},\\
J\textbf{M}&=&-c^{-1}\left((d_{ij}M_{j})_{,i}+b_{ji}M_{i,j}\right),\\
L_s \textbf{u}&=&\ell_{si} (B_{klji}u_{k,lj})_{,i},\\
Z_s\theta&=&-\ell_{si}\left( (b_{ji}\theta)_{,j}+d_{ij}\theta_{,j}\right),\\
N_s\textbf{R}&=&\ell_{si} (C_{jikl}R_{k,l})_{,j},\\
P_s\textbf{M}&=&\ell_{si} (\tilde{C}_{jikl}\theta_{k,l})_{,j},\\
\end{array}\ee  where $\ell_{si}$ is defined by $\ell_{si}c_{ij}=\delta_{sj}$.
We consider the matrix operator $\mathscr{A} $ on $\mathscr{H}$
 by \be \label {matrix2} \mathscr{A}= \left(
{\begin{array}{*{20}c}
  0& {\bf Id} &    0&  0&  0&  0  \\
  {\bf A}&  { 0}& 0&  {\bf B}&  {\bf C} & 0 \\
0&  { 0}& 0&{ Id} &    0&   0\\
0&  { D}& E&{ G}&   0&J   \\
0&  0&0&  0&0&   \textbf{Id} \\
      \textbf{L}& 0&0& \textbf{Z}&\textbf{N}&  \textbf{P}\\
\end{array}}\right)\ee
 where
$\textbf{Id}$ and $Id$ are  the identity operators in the respective
spaces, $\textbf{A}  =(A_i)$, $\textbf{B}  =(B_i)$, $\textbf{C}
=(C_i)$, $\textbf{L}=L_i$, $\textbf{Z}=Z_i$, $\textbf{N}=N_i$ and
$\textbf{P}=P_i$. The domain of $\mathscr{A}$ is \[\mathscr{D}
=\mathscr{D} (\mathscr{A} )=\left\{(u_i, v_i,\tau,\theta,R_i,M_i)\in
\mathscr{H};\ \mathscr{A}(u_i, v_i,\tau,\theta,R_i,M_i)\in
\mathscr{H}\right\}.\]

It is clear that
\[({\textbf{W}}_0^{1,2}(\Omega)\cap{\textbf{W}}^{2,2}(\Omega))\times{\textbf{W}}_0^{1,2}(\Omega)\times({{W}}_0^{1,2}(\Omega)
\cap{{W}}^{2,2}(\Omega))\times{{W}}_0^{1,2}(\Omega)
\times({\textbf{W}}_0^{1,2}(\Omega)\cap{\textbf{W}}^{2,2}(\Omega))\times{\textbf{W}}_0^{1,2}(\Omega)\]
is a subset of $\mathscr{D}$ which is dense in $\mathscr{H}$.

In the frame of  type II theory, we have that $ \tilde{K}_{ij} = 0$
 and $\tilde{C}_{ijkl} = 0$, so that $G =
0$ and $\textbf{P}=\textbf{0}$ in Eqs. (\ref{matrix1}) and
(\ref{matrix2}).

The initial boundary value problem (\ref{linear01}), (\ref{IC1}) and
(\ref{BC})
 can be transformed into the following Cauchy problem in
the Hilbert space $\mathscr{H}$,
 \be \label{system}
 \frac{d\mathcal{U}(t)}{dt}={\mathscr{A}}
\mathcal{U}(t)+{\mathcal{F}}(t),\ \ \
\mathcal{U}(0)=\mathcal{U}_0,\ee where
\[\mathcal{U}=(u_i,
v_i,\tau,\theta,R_i,M_i),\ \ \ \ {\mathcal{F}}=(0,\rho f_i,0, \rho
s,0,\rho Q_i),\ \ \ \ \mathcal{U}_0=(u^0_i,
v^0_i,\tau^0,\theta^0,R^0_i,M^0_i).\]

Now, we use the theory of semigroups of linear operators to obtain
the existence of solutions to the Cauchy problem (\ref{system}).
\begin{lemma}  The operator $\mathscr{A}$ satisfies the inequality
\[<{\mathscr{A}}\mathcal{U},\mathcal{U}> \le
0\] for every $\mathcal{U} \in \mathscr{D}(\mathscr{A})$, solution
to (\ref{system}).
\end{lemma}
{\bf Proof:} Let $\mathcal{U}=(u_i, v_i,\tau,\theta,R_i,M_i) \in
\mathscr{D}(\mathscr{A})$. Using the divergence theorem and the
boundary conditions, we have
 \begin{eqnarray}<{\mathscr{A}}\mathcal{U},\mathcal{U}>&=&
\int_\Omega\Big[{\mathcal{W}}\Big((u_i, \tau,R_i),(v_i,
\theta,M_i)\Big)-v_{i,j}t_{ji}-\Phi_i
\theta_{,i}-\Lambda_{ji}M_{i,j}
\Big]dv\cr&=&-\int_\Omega(\tilde{K}_{ij}\theta_{,i}\theta_{,j}+
\tilde{C}_{ijkl}M_{l,k}M_{i,j})dv\,.\cr \nonumber\end{eqnarray} The
thesis follows from our hypothesis (\ref{tilda}). \hfill$\Box$

In the context  of type II theory ($ \tilde{K}_{ij} =
\tilde{C}_{ijkl} = 0$), this lemma implies
$<{\mathscr{A}}\mathcal{U},\mathcal{U}>=0$, which means conservation
of the energy \be \label{E} E(t)=\half\int_\Omega\Big(\rho
{v_i}{v_i} +c \theta^2+c_{ij} M_i M_j + A_{ijkl}e_{kl}
e_{ij}+2B_{ijkl}e_{ij}R_{k,l}+K_{ij}
\tau_{,i}\tau_{,j}+C_{ijkl}R_{i,j}R_{k,l} \Big)dv.\ee This confirms
the result obtained  in \citep{IQ2}.

It is worth remarking that this quantity is also  conserved even if
we do not impose conditions $(\textrm{i})-(\textrm{iii})$.
\begin{lemma}  The operator $\mathscr{A}$ has the property that \[\hbox{Range}({I -\mathscr{A}})=\mathscr{H}.\]
\end{lemma}
{\bf Proof:} Let $\mathcal{U}^*=(u_i^*, v_i^*, \tau^*,\theta^*,
R_i^*,M_i^*)\in \mathscr{H}$. We must prove that the equation
\[\mathcal{U}-{\mathscr{A}}\mathcal{U} =\mathcal{U}^*\] has a
solution $\mathcal{U}=(u_i, v_i,\tau,\theta,R_i,M_i) \in
\mathscr{D}$. This equation leads to the system \be\label{system1}
\begin{array}{lll}
 \textbf{{u}}^*&=&\textbf{{u}}-\textbf{{v}},\cr \tau^*&=&\tau-\theta,\cr \textbf{R}^*&=&\textbf{R}-\textbf{M},\\
 \textbf{{v}}^*&=&\textbf{{v}}-(\textbf{{A u}}+{B} \theta+\textbf{C R}),\\
\theta^*&=&\theta-(D \textbf{{v}}+E \tau+G \theta+J \textbf{M}),\\
\textbf{M}^*
  &=&\textbf{M}-(\textbf{L u}+\textbf{Z} \theta+\textbf{N R}+\textbf{P M}).\\
\end{array}\ee Substituting the first three  equations in the others, we
obtain

\be\label{sys2}\mathscr{R}(\textbf{{u}},\tau,\textbf{R})=(
\textbf{{u}}^*+\textbf{{v}}^*-{B}\tau^*, \theta^*+\tau^*-D
{\textbf{u}}^*-G \tau^*- J\textbf{R}^*,
\textbf{R}^*+\textbf{M}^*-\textbf{Z} \tau^*- \textbf{P
}\textbf{R}^*)\ee where
\[\mathscr{R}= \left( {\begin{array}{*{20}c}
  {\bf Id-A} &    -B&  -\textbf{C}  \\
      {-{{D}}}&  Id-(E+G)& -J\\
-\textbf{{L}}&-\textbf{{Z}}&\textbf{Id}-(\textbf{{N}}+\textbf{P})\\
\end{array}}\right).\]

 To solve the system
\eqref{sys2} we introduce the following bilinear form on
${\textbf{W}}_0^{1,2}(\Omega)\times{{W}}_0^{1,2}(\Omega)\times
{\textbf{W}}_0^{1,2}(\Omega)$,
\[ {\mathscr{B}}[(\textbf{{u}},\tau,\textbf{R}),(\textbf{{u}}',\tau',\textbf{R}')]
=<\mathscr{R}(\textbf{{u}},\tau,\textbf{R}),(\textbf{{u}}',\tau',\textbf{R}')>_{{\bf{L}}^{2}\times
{{L}}^{2}\times{\bf{L}}^{2}}.\]

 A
direct calculation shows that ${\mathscr{B}}$ is bounded. Using the
divergence theorem, we have
\[ {\mathscr{B}}[(\textbf{{u}},\tau,\textbf{R}),(\textbf{{u}},\tau,\textbf{R})]
=\int_\Omega\Big[\rho u_i u_i+ A_{ijkl}e_{kl} e_{ij}+c
\tau^2+c_{ij}R_iR_j+\tilde{K}_{ij}\tau_{,i}\tau_{,j}+
\tilde{C}_{ijkl}R_{l,k}R_{i,j} +2\mathcal{W}[(\textbf{u},
\tau,\textbf{R}),(\textbf{u}, \tau,\textbf{R})]\Big]dv.\]
 In view of our assumptions on the constitutive
coefficients, we see that ${\mathscr{B}}$ is coercive on $
{\textbf{W}}^{-1,2}(\Omega)\times {{W}}^{-1,2}(\Omega)\times
{\textbf{W}}^{-1,2}(\Omega)$. On the other hand, it is easy to see
that the vector
\[(
\textbf{{u}}^*+\textbf{{v}}^*-{B}\tau^*, \theta^*+\tau^*-D
\textbf{{u}}^*-G \tau^*- J\textbf{R}^*,
\textbf{R}^*+\textbf{M}^*-\textbf{Z} \tau^*- \textbf{P
}\textbf{R}^*)\] lies in $ {\textbf{W}}^{-1,2}(\Omega)\times
{{W}}^{-1,2}(\Omega)\times {\textbf{{W}}}^{-1,2}(\Omega)$. Hence the
Lax-Milgram theorem \citep{[22]} implies the existence of
$(\textbf{u},\tau,\textbf{R})\in {\textbf{W}}_0^{1,2}(\Omega)\times
{{W}}_0^{1,2}(\Omega)\times {\textbf{{W}}}_0^{1,2}(\Omega)$ which
solves equation (\ref{sys2}). Now, we may also conclude the
existence of $\textbf{v}\in {\textbf{W}}_0^{1,2}(\Omega)$, $\theta
\in {{W}}_0^{1,2}(\Omega)$ and $\textbf{M} \in
{\textbf{{W}}}_0^{1,2}(\Omega)$ solving system
(\ref{system1}).\hfill$\Box$

The previous lemmas lead to next theorem.
\begin{theorem}
The operator $\mathscr{A} $ generates a semigroup of contraction in
$\mathscr{H}$.
\end{theorem}
{\bf Proof.} The proof follows from Lumer-Phillips corollary to the
Hille-Yosida theorem \citep{[23]}.\hfill$\Box$

 It is worth remarking that this
theorem implies that the dynamical system generated by the equations
of thermoelasticity with diffusion of type III  (or type II) are
stable in the sense of Lyapunov.
\begin{theorem}
Assume that the conditions $(\textrm{i})-(\textrm{iii})$ hold,
$f_i,\  s,\ Q_i \in C^1([0,\infty),L^2)\cap
C^0([0,\infty),{{W}}_0^{1,2})$ and $\mathcal{U}_0$ is in the domain
of the operator $\mathscr{A} $. Then, there exists a unique solution
$\mathcal{U}(t)\in C^1([0,\infty),\mathscr{H})\cap
C^0([0,\infty),\mathscr{D(A)})$ to the problem
(\ref{system}).\end{theorem} Since the solutions are defined by
means of a semigroup of contraction, we have the estimate
\[\| \mathcal{U}(t)\|_\mathscr{H}  \le \| \mathcal{U}_0 \|_\mathscr{H}+\int_0^t\Big(
\| {f_i}(\xi)\|_{\bf{L}^2}+\| s(\xi)\|_{{L}^2} +\|
 {Q_i}(\xi)\|_{\bf{L}^2}\Big)d\xi\] which proves the continuous dependence
of the solutions upon initial data and body loads. Thus, under
assumptions $(\textrm{i})-(\textrm{iii})$ the problem of linear
thermoelasticity
 with microtemperatures of type III (or type II)
is well posed.
\section{Asymptotic behavior of solutions}
In this section we study the asymptotic behavior of solutions, whose
existence has been  proved previously. We consider in the following
sections
 the homogenous case (${f_i}$=0,\ $s$=0,\ $Q_i=0$). In particular
we are interested in the relation between dissipation effects and
time decay of solutions. Therefore, we will  continue to assume that
the assumptions $(\textrm{i})-(\textrm{iii})$  considered in the
previous  section hold. However it is  worth noting that the results
for this section only hold for type III theory.

To this end, we recall that for a  semigroup of contraction,   the
precompact orbits tend to the $\omega-$limit  sets if its generator
$\mathscr{A}$ has only the fixed point ${\bf 0}$ (see \citep{[24]})
and the structure of the $\omega-$limit sets is determined by the
eigenvectors of eigenvalue $i \lambda$ (where $\lambda$ is a real
number) in the closed subspace
\[{\mathscr{L}}=\ll\{\mathcal{U}\in \mathscr{H};\
<\mathscr{A}\mathcal{U},\mathcal{U}>=0\}\gg\,,
\]
where $\ll \mathcal{C} \gg$ denotes the closed vectorial subspace
generated by the set $\mathcal{C}$.

From the assumptions $\textrm{(i)-(iii)}$ it is easy to check that
$\mathscr{A}^{-1}({\textbf{0}})={\textbf{0}}$, while the
precompactness of the orbits starting in $\mathscr{D}$ is a
consequence of the following Lemma  \citep{[23]}.
\begin{lemma}  The operator $({I- \mathscr{A}})^{-1}$ is compact.
\end{lemma}
{\bf Proof.} Let $(\hat{\textbf{u}}_n, \hat{\textbf{v}}_n,
 \hat{\tau}_n,\hat{\theta}_n, \hat{\textbf{R}}_n,\hat{\textbf{M}}_n)$ be a bounded sequence in
 $\mathscr{H}$ and let $U_n=({\textbf{u}}_n, {\textbf{v}}_n,
{\tau}_n, {\theta}_n,\textbf{{R}}_n, \textbf{{M}}_n)$ be the
sequence of the respective
 solutions to the system (\ref{system1}). We have
\be\label{PP}\mathscr{P}[(U_n,U_n)]=\mathscr{R}[({\textbf{u}}_n,
{\tau}_n, \textbf{{R}}_n),({\textbf{u}}_n, {\tau}_n,
\textbf{{R}}_n)]\leq \textrm{Constant}\times
\mathscr{P}[(U_n,U_n)]^{\frac{1}{2}}.\ee Inequality \eqref{PP}
implies that $({\textbf{u}}_n, {\tau}_n, \textbf{{R}}_n)$  is a
bounded sequence in $ {\textbf{W}}_0^{1,2}(\Omega)\times
{{W}}_0^{1,2}(\Omega)\times {\textbf{{W}}}_0^{1,2}(\Omega)$. The
theorem of Rellich-Kondrasov \citep{[25]} implies that there is
exists a subsequence converging in ${\bf{L}}^{2}(\Omega)\times
{{L}}^{2}(\Omega)\times {\textbf{{L}}}^{2}(\Omega)$. In a similar
way
\[ \textbf{v}_{n_j}=\textbf{u}_{n_j}-\hat{\textbf{u}}_{n_j},\ \ \
\theta_{n_j}={\tau}_{n_j}-\hat{\tau}_{n_j},\ \ \
\textbf{M}_{n_j}={\textbf{R}}_{n_j}-\hat{\textbf{R}}_{n_j}\] has a
sub-sequence converging in
${{{\bf{L}}^{2}(\Omega)\times{{L}}^{2}(\Omega)\times
{\textbf{{L}}}^{2}(\Omega)}}$. Thus we conclude the existence of a
sub-sequence
\[({\textbf{u}}_{n_{j_k}},{\textbf{v}}_{n_{j_k}},
{\tau}_{n_{j_k}},
{\theta}_{n_{j_k}},{\textbf{R}}_{n_{j_k}},\textbf{M}_{n_{j_k}})\]
which converges in $\mathscr{H}$. \hfill$\Box$

Now, we can state a theorem on the asymptotic behavior of solutions.
\begin{theorem}  Let $\mathcal{U}_0=({\bf{u}}^0,{ \bf{v}}^0,\tau^0,\theta^0,\textbf{R}^0,\textbf{M}^0)\in\mathscr{D(A)}$
and $\mathcal{U}(t)$ be the solution to the problem (\ref{system})
with $\mathcal{F}=0$.
 Then
\[ \tau(t) \to 0\ \hbox{as}\ t \to \infty\ \hbox{in}\ W_0^{1,2}(\Omega)\quad \hbox{and}
 \quad\theta(t) \to 0\ \hbox{as}\ t \to \infty\ \hbox{in}\ L^2(\Omega). \]
Moreover
\[ {\bf{u}}(t)  \to 0\ \hbox{as}\ t \to \infty\
\hbox{in}\ {\bf{W}_0^{1,2}}(\Omega) \quad \hbox{and}
 \quad{\bf{v}}(t)    \to 0\ \hbox{as}\ t \to \infty\ \hbox{in}\ {\bf{L}}^2(\Omega) \]
 \[ {\bf{R}}(t)  \to 0\ \hbox{as}\ t \to \infty\
\hbox{in}\ {\bf{W}_0^{1,2}}(\Omega) \quad \hbox{and}
 \quad{\bf{M}}(t)    \to 0\ \hbox{as}\ t \to \infty\ \hbox{in}\ {\bf{L}}^2(\Omega) \]
whenever the system
\begin{equation}
\label{subsystem}
\begin{array}{lll}
\bf{Au}+\lambda^2\bf{u}&=&0 \ \hbox{in}\ \Omega\\   D\bf{v}&=&0 \ \hbox{in}\ \Omega \\  \textbf{L}\bf{u}&=&0\ \hbox{in}\ \Omega\\
\bf{u}&=&0\ \hbox{on}\ \p \Omega
\end{array}
\end{equation}
 has only the null solution.\end{theorem}
{\bf Proof.} To prove the theorem we have to study the structure of
the $\omega-$limit set. Thus, we must study the equation \be
\label{hat} \hat{\mathscr{A}}\mathcal{U}=i \lambda \mathcal{U}\ee
for some real number $\lambda$, where $\mathcal{U}\in
\mathscr{D}(\hat{\mathscr{A}})$ and
$\hat{\mathscr{A}}={\mathscr{A}}_{|\mathscr{L}}$ is the generator of
a group on $\mathscr{L}$. If $\mathcal{U} \in \mathscr{L}$ then
$<\mathscr{A}\mathcal{U},\mathcal{U}>=0$. Under the assumption that
the tensors $\tilde{K}_{ij}$ and $\tilde{C}_{ijkl}$ are definite
positive, it follows that $\theta=M_i=0$ and then $\tau=R_i=0$.
Thus, the asymptotic behavior of the temperature and the
microtemperatures is proved. Now, Eq. (\ref{hat}) can be rewritten
as\[ \left( {\begin{array}{*{20}c}
  0& {\bf Id} &    0&  0&  0&  0  \\
  {\bf A}&  { 0}& 0&  {\bf B}&  {\bf C}& 0  \\
0&  { 0}& 0&{ Id} &    0&   0\\
0&  { D}& E&{ G}&   0&J   \\
0&  0&0&  0&0&   \textbf{Id} \\
      \textbf{L}& 0&0& \textbf{Z}&\textbf{N}&  \textbf{P}\\
\end{array}}\right)\left( {\begin{array}{*{20}c}
   {\bf u} &   \\
  {\bf v}&   \\
   0&\\
0&  \\
0&  \\
    0& \\
\end{array}}\right)=i\lambda\left( {\begin{array}{*{20}c}
   {\bf u} &   \\
  {\bf v}&   \\
   0&\\
0&  \\
0&  \\
    0& \\
\end{array}}\right).\]
Introducing the first equation into the others we obtain the system
(\ref{subsystem}).

If the system (\ref{subsystem}) has only the trivial solution, we
obtain $\omega-$limit$(\mathcal{U}_0) = 0$ and ${\bf{u}}(t)  \to 0\
\hbox{in}\ {\bf{W}}_0^{1,2}(\Omega) $ and $ {\bf{v}}(t)    \to 0\
\hbox{in}\ {\bf{L}}^2(\Omega) $ when $ t \to \infty$.\hfill$\Box$

\section{Impossibility of localization in time}In  previous sections
we have proved that the  solutions of type III theory are stable
asymptotically and also in the sense of Lyapunov. A natural question
is to ask if the decay is fast enough to guarantee that the
solutions vanishe in a finite time. In fact, when the dissipation
mechanism in a system is sufficiently strong, the localization of
solutions in the time variable can hold. This means that the decay
of the solutions is sufficiently fast to guarantee that they vanish
after a finite time.

In the context of Green-Naghdi thermoelasticity of type III,
\citet{[11]} has shown that the thermal dissipation is
not strong enough to obtain  the localization in time of the
solutions. In this section assume the quadratic form (\ref{tilda})
positive definite and prove that the further dissipation effects due
to the microtemperatures of type III are not sufficiently strong to
guarantee that the thermomechanical deformations vanish after a
finite interval of time. This means that, in absence of sources, the
only solution for the evolutive problem that vanishes after a finite
time is the null solution, that is the following theorem holds.
\begin{theorem} Let $(u_i, \tau,{R}_i)$ be a solution of the system (\ref{linear01}), (\ref{IC1}) and  (\ref{BC}) which vanishes after a finite time $t_0$. Then
$(u_i, \tau,{R}_i)\equiv (0,0,0)$ for every $t \ge 0$.
\end{theorem} In order to prove this theorem, generalizing the
technique used in \citep{[11]}, we
 show the uniqueness of solutions for the related  backward in time problem.
Backward in time problems are relevant from the mechanical point of
view when we want to have some information about what happened in
the past by means of the information that we have at this moment.

For our model, the system of equations which govern the backward in
time problem is given by
\begin{align}\begin{split}\label{back}
 \rho \ddot u_i &=(A_{jikl}e_{kl}+a_{ji}\dot \tau+B_{jikl}R_{l,k})_{,j}+\rho f_i,
\\
 c\ddot \tau&=a_{ij}\dot e_{ij}+(d_{ij}\dot R_{j}+K_{ij}\tau_{,j}-\tilde{K}_{ij}\dot\tau_{,j})_{,i} +b_{ij}\dot R_{i,j} +\rho s,
\\
c_{ij}\ddot R_{j}&= (B_{klji}e_{kl}+b_{ji}\dot
\tau+C_{jikl}R_{k,l}-\tilde{C}_{jikl}\dot R_{k,l})_{,j}+d_{ij}\dot
\tau_{,j}+\rho Q_i.
\end{split}\end{align}

\begin{proposition} [Uniqueness]   Let  $(u_i, \tau,{R}_i)$ be a solution
to the problem (\ref{back}), (\ref{BC}) with null initial data and
sources. Then $(u_i, \tau,{R}_i)=(0,0,0)$ for every $t\ge 0$.
\end{proposition}
{\bf Proof. } Let us introduce the following functionals
\begin{eqnarray*}E_1(t)&=&\half\int_\Omega\Big(\rho
\dot{u_i}\dot{u_i} +c \dot\tau^2+c_{ij} \dot R_i \dot R_j +
A_{ijkl}u_{i,j} u_{k,l}+2B_{ijkl}u_{i,j}R_{k,l}+K_{ij}
\tau_{,i}\tau_{,j}+C_{ijkl}R_{i,j}R_{k,l} \Big)dv\,,\\
 E_2(t)&=&\half\int_\Omega\Big(\rho
\dot{u_i}\dot{u_i} -c \dot\tau^2-c_{ij} \dot R_i \dot R_j +
A_{ijkl}u_{i,j} u_{k,l}-K_{ij}
\tau_{,i}\tau_{,j}-C_{ijkl}R_{i,j}R_{k,l}
\Big)dv\,,\\E_3(t)&=&\int_\Omega\Big(\rho {u_i}\dot{u_i}-c
\dot\tau\tau-c_{ij} \dot R_i R_j+\half \tilde{K}_{ij}
 \tau_{,i} \tau_{,j}+\half\tilde{C}_{ijkl}R_{i,j}R_{k,l}+a_{ij}\tau u_{i,j} \Big)dv\,,\\
 \end{eqnarray*}
and compute their time derivatives. By multiplying the  first
equation of (\ref{back}) by $ \dot u_i$, the second one by $\dot
\tau$  and the third one by $\dot R_i$, we get
\[\dot{E}_1(t)=\int_\Omega\Big(\tilde{K}_{ij}
\dot \tau_{,i}\dot \tau_{,j}+\tilde{C}_{ijkl}\dot   R_{i,j}\dot
R_{k,l} \Big)dv.\] If  we  multiply the  first equation of
(\ref{back}) by $\dot u_i$, the second one by $-\dot  \tau$ and the
third one by $-\dot R_i$, we have\begin{eqnarray*}
\dot{E}_2(t)&=&-\int_\Omega\Big(-2a_{ij}\dot u_{i}\dot\tau_{,j}+
B_{ijkl} (\dot u_{i,j}  R_{k,l}-u_{k,l}\dot R_{i,j})+\tilde{K}_{ij}
\dot \tau_{,i}\dot \tau_{,j}+\tilde{C}_{ijkl}\dot R_{k,l}\dot
R_{i,j}\Big)dv\end{eqnarray*} and, finally, if  we  multiply the
first equation of (\ref{back}) by $-  u_i$, the second one by $
\tau$  and the third one by $R_i$, we obtain$$
\dot{E}_3(t)=-\int_\Omega\Big( A_{ijkl}u_{k,l}  u_{i,j} +b_{ij}(\tau
\dot  R_{i,j}-\dot \tau R_{i,j})+d_{ij}(\dot \tau_{,j}
R_i-\tau_{,i}\dot R_j)+c \dot\tau^2+c_{ij} \dot R_i \dot R_j -K_{ij}
\tau_{,i}\tau_{,j}-{C}_{ijkl} R_{i,j} R_{k,l} -\rho\dot u_i\dot
u_i\Big)dv.$$

 Moreover,  a well-known identity
for type III thermoelasticity (see equation (3.9) in \citep{[10]})
for our model becomes
 \be \int_\Omega\Big(A_{ijkl} u_{i,j} u_{k,l}+c \dot\tau^2+c_{ij} \dot R_i \dot R_j\Big)dv
 = \int_\Omega\Big( \rho\dot u_i\dot u_i+K_{ij}
\tau_{,i}\tau_{,j}+{C}_{ijkl} R_{i,j} R_{k,l} \Big)dv.\ee
 Then we
have
\[{E}_2(t)=\int_\Omega\Big(
A_{ijkl}u_{i,j} u_{k,l}-K_{ij}
\tau_{,i}\tau_{,j}-C_{ijkl}R_{i,j}R_{k,l}\Big)dv
\]
and
\[\dot{E}_3(t)=-\int_\Omega\Big(b_{ij}(\tau \dot  R_{i,j}-\dot \tau
R_{i,j})+d_{ij}(\dot \tau_{,j} R_i-\tau_{,i}\dot R_j)\Big)dv.
\]
We consider the function
\begin{eqnarray*}  \mathcal{E}(t)&=&\int_0^t \left[\epsilon E_1(s)+ {E_2}(s)+\lambda E_3(s)\right]ds
\\&=&\frac12 \int_0^t \int_\Omega \left[\epsilon\rho \dot{u_i}\dot{u_i}+\epsilon c \dot\tau^2+\epsilon c_{ij}\dot   R_i \dot R_j
+(\epsilon+2)A_{ijkl} u_{i,j} u_{k,l} +2\epsilon
B_{ijkl}u_{i,j}R_{k,l}\right]dv\,ds
\\&+&\frac12 \int_0^t \int_\Omega\left\{[\lambda\tilde{K}_{ij}+(\epsilon-2){K}_{ij}] \tau_{,i} \tau_{,j}
+[\lambda\tilde{C}_{ijkl}+(\epsilon-2){C}_{ijkl}]  R_{i,j}
R_{k,l}\right\}dv\,ds
\\&+&\lambda \int_0^t \int_\Omega\left[\rho {u_i}\dot{u_i}-c \dot\tau\tau-c_{ij}  R_i \dot R_j+a_{ij}\tau u_{i,j} \right]dv\,ds
\end{eqnarray*}
where $\epsilon$ and $\lambda$ are positive suitable constants such
that the quadratic form
\[
\int_\Omega\left\{[\lambda\tilde{K}_{ij}+(\epsilon-2){K}_{ij}]
\tau_{,i} \tau_{,j}
+[\lambda\tilde{C}_{ijkl}+(\epsilon-2){C}_{ijkl}]  R_{i,j}
R_{k,l}\right\}\,dv
\]
is positive definite.
 By using the null initial data hypothesis and the Poincar\'e inequality we have
 \[\lambda \int_0^t \int_\Omega[\rho {u_i}\dot{u_i}-c \dot\tau\tau-c_{ij}  R_i \dot R_j]dv\,ds\leq
\frac{\epsilon}4 \int_0^t \int_\Omega[\rho \dot{u_i}\dot{u_i}+c
\dot\tau^2+ c_{ij} \dot  R_i \dot R_j ]\,dv\,ds
 \]
 for any $t\leq t_0$, where $t_0$ is a positive time which depends on $\lambda$, $\epsilon$ and the constitutive coefficients. Therefore $\mathcal{E}(t)$ is a positive definite quadratic form for $0\leq t\leq t_0$, in particular
 \begin{eqnarray}\label{calE} \mathcal{E}(t)&\geq&\frac14 \int_0^t \int_\Omega\left[\epsilon\rho \dot{u_i}\dot{u_i}
 +\epsilon c \dot\tau^2+\epsilon c_{ij}\dot   R_i \dot R_j+(\epsilon+2)A_{ijkl} u_{i,j} u_{k,l}
 +2\epsilon B_{ijkl}u_{i,j}R_{k,l}+2\lambda a_{ij}\tau
 u_{i,j}\right]dv\,ds\nonumber
\\&+&\frac14 \int_0^t \int_\Omega\left\{[\lambda\tilde{K}_{ij}+(\epsilon-2){K}_{ij}] \tau_{,i} \tau_{,j}
+[\lambda\tilde{C}_{ijkl}+(\epsilon-2){C}_{ijkl}]  R_{i,j}
R_{k,l}\right\}dv\,ds.
\end{eqnarray}
Moreover, recalling the null initial data assumption, we have
 \begin{eqnarray*} \dot{\mathcal{E}}(t)&=& (\epsilon-1)\int_0^t\int_\Omega\Big(\tilde{K}_{ij}
\dot \tau_{,i}\dot \tau_{,j}+\tilde{C}_{ijkl}\dot   R_{i,j}\dot
R_{k,l} \Big)\,dv\,ds-\int_0^t\int_\Omega\Big(-2a_{ij}\dot
u_{i}\dot\tau_{,j}+ B_{ijkl} (\dot u_{i,j} R_{k,l}-u_{k,l}\dot
R_{i,j})
\\&+&
\lambda b_{ij}(\tau \dot R_{i,j}-\dot \tau R_{i,j})+\lambda
d_{ij}(\dot \tau_{,j} R_i-\tau_{,i}\dot R_j)\Big)dvds.
 \end{eqnarray*}
Choosing $0<\epsilon<1$ and using the inequality of arithmetic and
geometric means, we have
 \begin{eqnarray*}
&& \left|\int_0^t\int_\Omega\Big(-2 a_{ij}\dot u_{i}\dot\tau_{,j}+
B_{ijkl} (\dot u_{i,j} R_{k,l}-u_{k,l}\dot R_{i,j}) + \lambda
b_{ij}(\tau \dot R_{i,j}-\dot \tau R_{i,j})+\lambda d_{ij}(\dot
\tau_{,j} R_i-\tau_{,i}\dot R_j)\Big)dvds\right|\leq
\\
&&(1-\epsilon)\int_0^t\int_\Omega\Big(\tilde{K}_{ij} \dot
\tau_{,i}\dot \tau_{,j}+\tilde{C}_{ijkl}\dot   R_{i,j}\dot R_{k,l}
\Big)\,dv\,ds+K_1  \int_0^t \int_\Omega \rho
\dot{u_i}\dot{u_i}\,dv\,ds +K_2 \int_0^t \int_\Omega[c
\dot\tau^2+c_{ij} \dot R_i \dot R_j]dv\,ds
\\&&+K_3 \int_0^t \int_\Omega A_{ijkl} u_{i,j} u_{k,l} dv\,ds+K_4
\int_0^t \int_\Omega B_{ijkl}u_{i,j}R_{k,l} dv\,ds+K_5 \int_0^t
\int_\Omega a_{ij}u_{i,j}\tau dv\,ds
\\
&&+\frac12
\int_0^t\int_\Omega\left\{[\lambda\tilde{K}_{ij}+(\epsilon-2){K}_{ij}]
\tau_{,i} \tau_{,j}
+[\lambda\tilde{C}_{ijkl}+(\epsilon-2){C}_{ijkl}]  R_{i,j}
R_{k,l}\right\}\,dv \,ds
\end{eqnarray*}
where the positive constants $K_i$ can be calculated by standard
methods, so that
 \begin{eqnarray}   \label{dotcalE} \dot{ \mathcal{E}}(t)&\leq&K\int_0^t \int_\Omega \left\{\epsilon\rho \dot{u_i}\dot{u_i}+
 \epsilon c \dot\tau^2+\epsilon c_{ij} \dot R_i \dot
R_j+(\epsilon+2)A_{ijkl} u_{i,j} u_{k,l}+2\epsilon
B_{ijkl}u_{i,j}R_{k,l}+2\lambda a_{ij}u_{i,j}\tau
\right\}dv\,ds\nonumber
\\&+&K \int_0^t \int_\Omega\left\{[\lambda\tilde{K}_{ij}+(\epsilon-2){K}_{ij}] \tau_{,i} \tau_{,j}
+[\lambda\tilde{C}_{ijkl}+(\epsilon-2){C}_{ijkl}]  R_{i,j}
R_{k,l}\right\}dv\,ds
\end{eqnarray}
with $K=\max\{\frac12,K_1,K_2,K_3,K_4,K_5\}$. Inequalities
(\ref{calE}) and (\ref{dotcalE}) yield
\[ \dot{\mathcal{E}}(t)\le 4K \mathcal{E}(t)\,,\qquad 0\le t \le t_0\,.\]
This inequality and the null initial data imply
$\mathcal{E}(t)\equiv 0$ if  $0\le t \le t_0$. Reiterating this
argument on each subinterval $[(n-1)t_0 , nt_0]$ we obtain
$\mathcal{E}(t)\equiv 0$ for  $ t \geq0$.

If we take into account the definition of $\mathcal{E}(t)$,
the uniqueness result is proved. \hfill$\Box$

\section*{Acknowledgment}
Part of this work was done when the first author visited the
Dipartimento di Ingegneria Industriale, Universit\`a degli Studi di
Salerno, in June 2016 and February 2017. He thanks their hospitality
and financial support  through FARB  2015.

\section*{References}

\bibliography{references}

\end{document}